\documentclass[12pt,reqno]{amsart}
\usepackage[utf8]{inputenc}
\usepackage{amsmath,amsthm,amssymb}
\usepackage{mathtext}
\usepackage[T1,T2A]{fontenc}
\usepackage[english]{babel}
\usepackage{tikz}

\newcommand{\KK}{\mathbb K}

\newcommand{\ZZ}{\mathbb Z}

\newcommand{\Ker}{\mathrm{Ker}\,}

\newcommand{\LND}{\mathrm{LND}}
\renewcommand{\phi}{\varphi}

\usepackage{graphicx}
\usepackage[dvipsnames]{xcolor}
\usepackage[colorlinks=true, allcolors=red]{hyperref}
\usepackage{cite}

\theoremstyle{plain}
\newtheorem{theorem}{Theorem}
\newtheorem{lemma}{Lemma}

\theoremstyle{remark}
\newtheorem{remark}{Remark}
\theoremstyle{definition}

\newtheorem{question}{Question}

\title{}

\author{Veronika Kikteva}
\address{HSE University, Faculty of Computer Science, 11 Pokrovsky Bulvar, Moscow, 109028, Russia}
\email{VVKikteva@yandex.ru}

\thanks{The work was supported by the Foundation for the Advancement of Theoretical Physics and Mathematics “BASIS”}



\begin{document}
\title{Rigidity does not imply substitution rigidity}
\maketitle
\begin{abstract}
We obtain a negative answer to a question of Freudenburg concerning rigidity conditions for a polynomial. We construct a polynomial that does not lie in the kernel of any nonzero locally nilpotent derivation on the polynomial ring, but there exist an affine domain, algebraically independent elements and a locally nilpotent derivation such that the result of substituting elements into the polynomial lies in the kernel of that locally nilpotent derivation, but at least one of the substituted elements does not lie in that kernel.
\end{abstract}

\section{Introduction}
Let $\mathbb{K}$ be an algebraically closed field of characteristic zero. All rings considered in this paper are associative, commutative, and have an identity element. Let $A$ be an algebra over the field $\mathbb K$. A $\KK$-\textit{derivation} of $A$ is a $\KK$-linear map $D$: $A\to A$ satisfying the Leibniz rule, that is, for any elements $a,b\in A$, we have $D(ab)=aD(b)+bD(a)$. From now on all derivations are assumed to be $\KK$-derivations, unless otherwise is specified. A derivation $D$ of $A$ is called a \textit{locally nilpotent derivation} (\textit{LND} for short) if for every~$a\in A$ there exists a natural number $n$ such that $D ^n (a)=0$. We denote the set of locally nilpotent derivations of the algebra~$A$ by $\mathrm{LND}(A)$.

We say that a polynomial $P\in\mathbb{K}[X_1,\dots,X_m]$ is \textit{substitution-rigid} if for any affine $\mathbb{K}$-domain $B$, any algebraically independent elements $a_1,\dots,a_m\in B$, and any $D\in \mathrm{LND}(B)$, the condition $P(a_1,\dots,a_m)\in \Ker\,D$ implies
${a_1,\dots,a_m\in \Ker D}$.

We give several examples illustrating substitution rigidity. The kernel of an LND on a domain is factorially closed, that is, for any $\mathbb K$-domain $B$, any LND $D$ on $B$, and any nonzero elements $a,b\in B$ the condition $ab\in \Ker D$ implies $a,b\in\Ker D$; see~\cite[Principle 1 (a)]{F}. Hence, the monomials $X_1^{\alpha_1}\dots X_m^{\alpha_m}$, where $\alpha_i\geq 1$ for each $i$, are substitution-rigid. Another class of substitution-rigid polynomials consists of polynomials of the form ${c_1X_1^{\alpha_1}+c_2X_2^{\alpha_2}}$, where ${\alpha_1,\alpha_2\geq 2}$, ${c_1,c_2\in\mathbb{K}^{\times}}$ and $m=2$ by the Second AB Theorem, see~\cite[Theorem 2.50]{F}. A generalization of this result to $k$-nomials of the form
$$c_1 X_{11}^{\alpha_{11}}\dots X_{1n_1}^{\alpha_{1n_1}}+\dots +c_k X_{k1}^{\alpha_{k1}}\dots X_{kn_k}^{\alpha_{kn_k}},$$
where $k\geq 3$, $c_i\in\KK^{\times}$, $n_i\geq 1$, $\alpha_{ij}\geq 1$ and $$\sum_{i=1}^{k}\sum_{j=1}^{n_i} \frac{1}{\alpha_{ij}}\leq \frac{1}{k-1},$$ in the case $m=n_1+\dots +n_k$, was obtained in~\cite[Theorem 3.3]{K23}. Coordinates, i.e., polynomials that can be included in some coordinate system, are not substitution-rigid in the case ${m>1}$.

We say that a polynomial $P\in \mathbb{K}[X_1,\dots,X_m]$ is \textit{rigid} if it does not belong to the kernel of any nonzero $D\in\mathrm{LND}(\mathbb{K}[X_1,\dots, X_m])$. Note that this condition is not equivalent to the rigidity of the quotient algebra $\KK[X_1,\dots,X_m]/(P)$; see~\cite[Example 1]{K26}.

\begin{remark}
Rigidity follows from substitution rigidity. Suppose $P\in\mathbb{K}[X_1,\dots,X_m]$ is substitution-rigid. Take $B=\mathbb{K}[X_1,\dots,X_m]$. The elements $a_i:=X_i$ are algebraically independent. If $P$ belongs to the kernel of ${D\in\mathrm{LND}(\mathbb{K}[X_1,\dots,X_m])}$, then, since $P$ is substitution-rigid, we obtain $X_i\in\Ker D$, and hence $D$ is the zero derivation.
\end{remark}

The following question was posed in the book~\cite{F}.
\begin{question}(\cite[Question 11.9]{F})
\label{q}
Does substitution rigidity follow from rigidity?
\end{question}

In the case $m=2$ the answer is positive for irreducible polynomials. Kaliman and Makar-Limanov proved in~\cite[Corollary 2.1]{KML} that any irreducible polynomial $P\in\mathbb{C}[x,y]$ that is not a coordinate is substitution-rigid. For an irreducible polynomial in the case $m=2$, the condition of not being a coordinate is equivalent to rigidity; see~\cite[Section 11.7]{F}. In~\cite{K26} the equivalence of substitution rigidity and rigidity was proved for binomials.

In the present paper we give a negative answer to Question~\ref{q}.

\begin{theorem}\label{thm}
    For each $m\geq 3$ there is an irreducible polynomial $P\in\KK[X_1,\dots,X_m]$ such that $P$ is rigid, but not substitution-rigid.
\end{theorem}

 That is, $P$ does not lie in the kernel of any nonzero LND on $\KK[X_1,\dots,X_m]$, but there exist an affine domain $B$, algebraically independent elements $a_1,\dots,a_m\in B$ and ${D\in\LND(B)}$ such that $P(a_1, \dots,a_m)\in\Ker D$ and $a_1\not\in\Ker D$. Moreover, one can take $B=\KK[x_1,\dots,x_m]\simeq \KK^{[m]}$ and $D=\frac{\partial}{\partial x_3}$.

\textbf{Declaration of AI assistance.} The author used OpenAI’s GPT-6 Astra model to generate the example in the case $m=3$. The author independently checked all AI-assisted arguments and takes full responsibility for the correctness and originality of the manuscript.

The author is grateful to Sergey Gaifullin for constant
attention to this work.

\section{Example for $m=3$}
\label{sectm=3}
Consider a polynomial
\begin{equation}
\label{eqm=3}
    P=X^2+Y^3+Z^3\in\KK[X,Y,Z].
\end{equation}
It can easily be checked that $P$ is irreducible. By~\cite[Theorem 9.7(a)]{F}, the quotient algebra $\KK[X,Y,Z]/(P)$ is rigid, i.e., there are no nonzero LNDs on it. By~\cite[Remark 2]{K26}, the polynomial $P$ is rigid.

\begin{lemma}\label{lemm=3notsr}
The polynomial $P$ is not substitution-rigid.
\end{lemma}
\begin{proof}
Consider $B=\KK[x,y,z]\simeq \KK^{[3]}$ and $D=\frac{\partial}{\partial z}$. Put 
\begin{equation}
\label{eqtuvw}
t=y^3-x^3,\ u=x^3+y^3,\ v=tu=y^6-x^6,\ w=-vz^2-3xy=(x^6-y^6)z^2-3xy
\end{equation}
and define
\begin{equation}
\label{eqa123}
    a_1=vzw,\ a_2=xt+y^2w,\ a_3=yt-x^2w.
\end{equation}

The elements $a_1,a_2,a_3$ are algebraically independent. Indeed, consider their Jacobian matrix at $z=0$: $$\left(\frac{\partial(a_1,a_2,a_3)}{\partial(x,y,z)}\right)\Bigg|_{z=0}=\begin{pmatrix}
    0 & 0 & 3xy(x^6-y^6)\\
    - 4x^3 - 2y^3&- 6xy^2 & 0\\
    6x^2y & 2x^3 + 4y^3 &0
\end{pmatrix}.$$ Its determinant equals $-24xy(x^6-y^6)(x^3-y^3)^2\neq 0$. By~\cite[Lemma 3.10 (c)]{F}, the elements $a_1,a_2,a_3$ are algebraically independent.

Note that $$a_2^3+a_3^3=(xt+y^2w)^3+(yt-x^2w)^3=t^3(x^3+y^3)+3xytw^2(y^3+x^3)+w^3(y^6-x^6)=$$
$$=t^3u+3xyvw^2+w^3v=t^3u+vw^2(3xy+w)=t^3u-v^2w^2z^2=t^3u-a_1^2.$$

Therefore, $$a_1^2+a_2^3+a_3^3=t^3u=(y^3-x^3)^3(x^3+y^3)\in \Ker  D.$$

However, $$D(a_1)=v(zD(w)+w)=v(z(-2vz)+w)= $$ $$ =v(-2vz^2-vz^2-3xy)=-3v(vz^2+xy)\neq 0.$$

Thus $P(a_1,a_2,a_3)\in\Ker D$, but $a_1\not\in\Ker D$. So,~$P$ is not substitution-rigid.
\end{proof}

\section{Generalization to the case $m\geq 4$}
\label{sectm>=4}
Define the polynomial 
\begin{equation}
\label{eqm=4}
P_4:=(X_1^2+X_2^3+X_3^3)^2+X_4^3\in\KK[X_1,\dots,X_4]
\end{equation}
and for each $m\geq 5$ put 
\begin{equation}
\label{eqm>=5}
P_m:=(X_1^2+X_2^3+X_3^3)^N+X_4^N+\dots+X_m^N\in\KK[X_1,\dots,X_m],
\end{equation}
where $N:=(m-2)(m-3)$.

\begin{lemma}
\label{lemm>=4irr}
    The polynomial $P_m$ is irreducible for any $m\geq 4$.
\end{lemma}
\begin{proof}
    \textbf{Case $m=4$.} Consider the following $\ZZ$-grading on $\KK[X_1,\dots,X_4]$: $\deg(X_1)=3$, $\deg(X_2)=\deg(X_3)=2$, $\deg(X_4)=4$. The polynomial $P_4$ is homogeneous of degree~12. Suppose $P_4=FG$ for some nonconstant polynomials $F$ and $G$. Then $F$ and $G$ are homogeneous and $$F(X_1,X_2,X_3,0)G(X_1,X_2,X_3,0)=(X_1^2+X_2^3+X_3^3)^2.$$ Since $X_1^2+X_2^3+X_3^3$ is irreducible, we conclude $$F(X_1,X_2,X_3,0)=c(X_1^2+X_2^3+X_3^3),$$ $$G(X_1,X_2,X_3,0)=c^{-1}(X_1^2+X_2^3+X_3^3)$$ for some $c\in\KK^{\times}$. Hence $$\deg(F)=\deg(G)=\deg(X_1^2+X_2^3+X_3^3)=6$$ and the degrees of $F$ and $G$ with respect to $X_4$ are at most 1. This implies that $$\deg_{X_4}(P_4)=3=\deg_{X_4}(F)+\deg_{X_4}(G)\leq 2,$$ a contradiction.

    \textbf{Case $m\geq 5$.} Consider $$H:=(X_1^2+X_2^3+X_3^3)^N+X_4^N+\dots+X_{m-1}^N,$$ so $P_m=X_m^N+H$. Let $Q$ be an irreducible divisor of $H$. We apply Eisenstein's criterion for~$P_m$, considered as a polynomial in $X_m$ over $\KK[X_1,\dots,X_{m-1}]$; see~\cite[Chapter IV, Theorem 3.1]{L}. Trivially, $Q$ divides $H$ and does not divide $1$. Note also that $Q^2$ does not divide~$H$, since $$\gcd\left(H,\frac{\partial H}{\partial X_4}\right)=\gcd((X_1^2+X_2^3+X_3^3)^N+X_4^N+\dots+X_{m-1}^N, NX_4^{N-1})=1$$ and $H$ is square-free. By Eisenstein's criterion the polynomial $P_m$ is irreducible.
\end{proof}

\begin{lemma}
\label{lemm>=4r}
    The polynomial $P_m$ is rigid for any $m\geq 4$.
\end{lemma}
\begin{proof}
Suppose $P_m$ is in the kernel of a nonzero $D\in\LND(\KK[X_1,\dots,X_m])$. If $m=4$, the elements $X_1^2+X_2^3+X_3^3$ and $X_4$ are in the kernel of $D$ by~\cite[Theorem 2.50 (a)]{F}. If $m\geq 5$ we apply~\cite[Theorem 3.3]{K23} and conclude that $X_1^2+X_2^3+X_3^3, X_4,\dots,X_m$ are in the kernel of $D$ as well, since $$(m-2)\frac{1}{(m-2)(m-3)}=\frac{1}{m-3}= \frac{1}{(m-2)-1}$$ and no nonempty subsum vanishes.

Extend the field $\KK$ by adjoining $X_4,\dots,X_m$, and consider its algebraic closure
$$\mathbb{L}:= \overline{\mathbb{K}(X_4,\dots,X_m)}.$$
The algebra $\KK[X_1,\dots,X_m]$ can be embedded into
$B'=\mathbb{L}[X_1,X_2,X_3]$
and $D$ extends to a nontrivial locally nilpotent $\mathbb{L}$-derivation $D'$ on $B'$ such that $D'(X_1^2+X_2^3+X_3^3)=0$. Since $X_1^2+X_2^3+X_3^3$ is rigid in the case $m=3$, we conclude that $D'=0$. Therefore, $D$ vanishes on $X_1,X_2,X_3$. Hence $X_1,\dots,X_m\in\Ker D$ and~$D$ is the zero derivation, a contradiction. This implies that $P_m$ is rigid.
\end{proof}

\begin{lemma}
\label{lemm>=4nsr}
The polynomial $P_m$ is not substitution-rigid.
\end{lemma}
\begin{proof}
Consider $B=\KK[x,y,z,a_4,\dots,a_m]\simeq \KK^{[m]}$ and an LND $D=\frac{\partial}{\partial z}$. Define elements $a_1,a_2,a_3$ by~(\ref{eqtuvw})--(\ref{eqa123}). The elements $a_1,\dots,a_m$ are algebraically independent, since their Jacobian matrix is block diagonal with non-degenerate blocks $\left(\frac{\partial(a_1,a_2,a_3)}{\partial(x,y,z)}\right)$ and the ${(m-3)\times (m-3)}$ identity matrix. As before $P_m(a_1,\dots,a_m)\in\Ker D$ and $a_1\not\in\Ker D$. Therefore, $P_m$ is not substitution-rigid.
\end{proof}

\begin{proof}[Proof of Theorem~\ref{thm}]
    For $m=3$ take the polynomial~(\ref{eqm=3}). It is irreducible and rigid as was shown in Section~\ref{sectm=3}. By Lemma~\ref{lemm=3notsr} $P$ is not substitution-rigid.

    For $m=4$ consider $P_4$ from~(\ref{eqm=4}) and for $m\geq 5$ take the polynomial defined by~(\ref{eqm>=5}). By Lemma~\ref{lemm>=4irr} $P_m$ is irreducible. By Lemmas~\ref{lemm>=4r} and~\ref{lemm>=4nsr} it is rigid, but not substitution-rigid.
    
\end{proof}


\begin{thebibliography}{99}
\bibitem{F}
Gene Freudenburg. \textit{Algebraic theory of locally nilpotent derivations}. Encyclopaedia Math. Sci.~136, Springer-Verlag, Berlin, Heidelberg, 2017

\bibitem{KML}
Shulim Kaliman and Leonid Makar-Limanov.
\textit{AK-invariant of affine domains}. Affine Algebraic Geometry. Osaka University Press. Osaka, Japan (2007), 231–-255

\bibitem{K23}
Veronika Kikteva. \textit{Generalization of the ABC Theorem on Locally Nilpotent Derivations}. Sib. Math. J. 64:5 (2023), 1167--1178

\bibitem{K26}
Veronika Kikteva. \textit{Rigidity conditions for binomials}. arXiv:2609.36135

\bibitem{L}
Serge Lang. \textit{Algebra}. Grad. Texts in Math.~211, Springer-Verlag, New York, 2002

\end{thebibliography}
\end{document}